\documentclass[12pt]{amsart}

\allowdisplaybreaks

\usepackage{fullpage}

\usepackage{amsfonts,amssymb,amscd,amsmath,enumerate,verbatim,color,xcolor}
\usepackage[latin1]{inputenc}
\usepackage{amscd}
\usepackage{latexsym}
\usepackage{tikz}

\usepackage{graphicx} 

\usepackage{mathptmx}

\usepackage{hyperref}
\hypersetup{colorlinks,linkcolor=blue ,citecolor=blue, urlcolor=blue}
\def\NZQ{\mathbb}               
\def\NN{{\NZQ N}}

\def\ZZ{{\NZQ Z}}
\def\RR{{\NZQ R}}

\def\Ghat{{\widehat{G}}}

\def\frk{\mathfrak}               

\def\Phi{{\frk N}}
\def\ab{{\mathbf a}}

\def\eb{{\mathbf e}}

\def\xb{{\mathbf x}}
\def\yb{{\mathbf y}}

\def\opn#1#2{\def#1{\operatorname{#2}}} 
\opn\gr{gr}

\def\Pc{{\mathcal P}}

\newtheorem{Theorem}{Theorem}[section]
\newtheorem{Lemma}[Theorem]{Lemma}
\newtheorem{Corollary}[Theorem]{Corollary}
\newtheorem{Proposition}[Theorem]{Proposition}

\theoremstyle{definition}
\newtheorem{Remark}[Theorem]{Remark}

\newtheorem{Example}[Theorem]{Example}

\newtheorem{Definition}[Theorem]{Definition}

\newtheorem{Conjecture}[Theorem]{Conjecture}

\let\epsilon\mu
\let\phi=\varphi
\let\kappa=\varkappa
\opn\dis{dis}
\opn\height{height}
\opn\dist{dist}
\def\pnt{{\raise0.5mm\hbox{\large\bf.}}}

\opn\Lex{Lex}
\opn\conv{conv}

\begin{document}

\title{Number of facets of symmetric edge polytopes arising from\\
join graphs}

\author{Aki Mori, Kenta Mori and Hidefumi Ohsugi}

\address{Aki Mori,
	Learning center, Institute for general education, Setsunan University, Neyagawa, Osaka, 572-8508, Japan} 
\email{aki.mori@setsunan.ac.jp}

\address{Kenta Mori,
	Department of Mathematical Sciences,
	School of Science,
	Kwansei Gakuin University,
	Sanda, Hyogo 669-1330, Japan}
\email{k-mori@kwansei.ac.jp}

\address{Hidefumi Ohsugi,
	Department of Mathematical Sciences,
	School of Science,
	Kwansei Gakuin University,
	Sanda, Hyogo 669-1330, Japan} 
\email{ohsugi@kwansei.ac.jp}


\begin{abstract}
Symmetric edge polytopes of graphs 
are important object in Ehrhart theory,
and have an application to Kuramoto models.
In the present paper, we study the upper and lower bounds
for the number of facets of
symmetric edge polytopes of connected graphs
conjectured by Braun and Bruegge.
In particular, we show that their conjecture is true
for any graph that is the join of two graphs (equivalently, 
for any connected graph whose complement graph is not connected).
It is known that any symmetric edge polytope is a centrally symmetric
reflexive polytope.
Hence our results give a partial answer to Nill's conjecture: the number of facets of a $d$-dimensional reflexive polytope is at most $6^{d/2}$.
\end{abstract}

\maketitle

\section{Introduction}

A {\it lattice polytope} $\Pc\subset \RR^d$ is a convex polytope 
all of whose vertices belong to $\ZZ^d$.
A $d$-dimensional lattice polytope $\Pc \subset \RR^d$ is called \textit{reflexive} if the origin of $\RR^d$ belongs to the interior of $\Pc$ and its dual polytope 
\[\Pc^\vee:=\{\yb \in \RR^d  :  \langle \xb,\yb \rangle \leq 1 \ \text{for all}\  \xb \in \Pc \}\]
is also a lattice polytope, where $\langle \xb,\yb \rangle$ is the usual inner product of $\RR^d$.
In general, we say that a lattice polytope is reflexive if it is unimodularly equivalent to a reflexive polytope.
It is known \cite{mirror}
that reflexive polytopes correspond to Gorenstein toric Fano varieties, and they are related
to mirror symmetry.
Let $N(\Pc)$ be the number of facets of a lattice polytope $\Pc$.
If $\Pc$ is reflexive, then $N(\Pc)$ is the number of vertices of the reflexive polytope $\Pc^\vee$.
The number $N(\Pc)$ is important when $\Pc$ is a $d$-dimensional reflexive polytope
since $N(\Pc) - (d+1)$ is the rank of the class group of the associated toric variety. 
Nill conjectured (a dual version of) the following.

\begin{Conjecture}[{\cite[Conjecture 5.2]{NillConj}}] \label{NConj}
Let $\Pc$ be a $d$-dimensional reflexive polytope.
Then $N(\Pc) \le 6^{d/2}$.
\end{Conjecture}

Nill \cite{Nillps} showed that Conjecture~\ref{NConj} is true for any
pseudo-symmetric reflexive simplicial $d$-dimensional polytope
and the maximum $6^{d/2}$ is attained if and only if 
$\Pc$ is a free sum of $d/2$ copies of del Pezzo polygons.

On the other hand, Higashitani \cite{Higashi} showed that
centrally symmetric simplicial reflexive polytopes are precisely the ``symmetric edge polytopes" of graphs without even cycles.
The definition of symmetric edge polytopes is as follows.
Let $G$ be a finite simple graph on the vertex set $[n]:=\{1,\ldots,n\}$ with the edge set $E(G)$. 
The \textit{symmetric edge polytope $\Pc_G$} of $G$ is the convex hull of
$
 \{ \pm (\eb_i - \eb_j) : \{i,j\} \in E(G) \},
$
where $\eb_i$ is the $i$-th unit coordinate vector in $\RR^n$. 
It is known that the symmetric edge polytope 
of a connected graph with $n$ vertices is a centrally symmetric
reflexive $(n-1)$-dimensional polytope.
Symmetric edge polytopes are studied in several different areas.
\begin{itemize}
\item[(a)]
Ehrhart theory:
The name ``symmetric edge polytope" was given in \cite{MHNOH}
in the study of Ehrhart theory.
Given a lattice polytope $\Pc \subset \RR^d$,
the {\it Ehrhart polynomial} of $\Pc$ is defined by
$E_\Pc(n) = | n \Pc \cap \ZZ^d|$ for $n \in \NN$
which is a polynomial in the variable $n$ of degree $\dim \Pc$.
It is known that the coefficients of the {\it $h^*$-polynomial}
$h^*_\Pc(\lambda)$
defined by $1+ \sum_{n=1}^\infty E_\Pc (n) \lambda^n = h^*_\Pc(\lambda)/(1-\lambda)^{d+1}$
are nonnegative integers.
Moreover it is known \cite{hibi} that a $d$-dimensional lattice polytope $\Pc$ is reflexive if and only if
$h^*_\Pc(\lambda)$ is palindromic, i.e., $h^*_\Pc (\lambda) = \lambda^d h^*_\Pc (1/\lambda)$.
One of the most important problems on palindromic $h^*$-polynomials is their
real-rootedness and gamma positivity.
In \cite{HJM}, the $h^*$-polynomial of the symmetric edge polytope $\Pc_G$ of the complete bipartite graph was given explicitly.
The facet description of the symmetric edge polytopes played an important role for the proof.
The gamma positivity of $h^*$-polynomial of symmetric edge polytopes
is studied in \cite{DJKV, KT, OT1, OT2}.
 D'Al\`i et al.~\cite{DJK},  
 gave a generalization of symmetric edge polytope to regular matroids, and showed that 
 two symmetric edge polytopes are unimodularly
 equivalent if and only if they correspond to the same graphic matroid.

\item[(b)]
Application to Kuramoto models:
Symmetric edge polytopes are known as {\it adjacency polytopes} (\cite{CDM}) which have
an application to Kuramoto models.
The normalized volume of the symmetric edge polytope is
an upper bound of the number of possible solutions in the Kuramoto equations.
In \cite{CDM, DDM}, explicit formulas of the normalized volumes of the symmetric
edge polytopes of certain classes of graphs are given by using
the facet descriptions of the symmetric edge polytopes.

\end{itemize}


In both (a) and (b), the facet descriptions of the symmetric edge polytopes
play important roles.
Motivated by its increasing importance,
Chen et al.~\cite{CDK} gave descriptions of the correspondence between faces of a symmetric edge polytope and face subgraphs
 of the underlying connected simple graph.
On the other hand, Braun and Bruegge \cite{BrBr, BBK} studied
upper and lower bounds for the number of the facets of symmetric edge polytopes.
Let $G_1$ and $G_2$ be graphs with exactly one common vertex.
Then the {\it $1$-sum} (called \textit{wedge} in \cite{BrBr}) 
of $G_1$ and $G_2$ is the union of $G_1$ and $G_2$.
The $1$-sum of several graphs are defined by a sequence of $1$-sums.
It is known \cite{BrBr} that $N(\Pc_G) = N(\Pc_{G_1})N(\Pc_{G_2})$ if $G$ is the $1$-sum of $G_1$ and $G_2$.
Let $K_n$ denote the complete graph with $n$ vertices,
and let $K_{\ell_1,\dots, \ell_s}$ denote the complete multipartite graph
on the vertex set $V_1 \sqcup \dots \sqcup V_s$ with $|V_i|=\ell_i$.
It is known \cite{HJM} that 
\begin{eqnarray}
    N(\Pc_{K_{\ell,m}}) &=& 2^\ell + 2^m-2,\label{2bu graph}\\
    N(\Pc_{K_{\ell_1,\dots,\ell_s}}) &= & 2^{\sum_{i=1}^s \ell_i} -  \sum_{i=1}^s ( 2^{\ell_i} -2) -2 \ \ \  \mbox{ if } s \ge 3.\label{tabu graph}
\end{eqnarray}
In particular, we have $N(\Pc_{K_n}) = 2^n-2$.
Braun and Bruegge \cite{BrBr} conjectured the following, and studied $N(\Pc_G)$ for sparse graphs $G$.
(Note that $2^{\frac{n}{2}+1}-2<3  \cdot 2^{\frac{n-1}{2}}-2$
and
$14 \cdot 6^{\frac{n}{2}-2} < 6^{\frac{n-1}{2}}$ for any $n \in \mathbb{N}$.)

\begin{Conjecture}[{\cite[Conjecture 2]{BrBr}}] \label{yoso}
Let $G$ be a connected graph with $n\ge 3$ vertices.
\begin{itemize}
    \item[(1)]
    If $n$ is odd, then we have
$3  \cdot 2^{\frac{n-1}{2}}-2 \le N(\Pc_G) \le 6^{\frac{n-1}{2}}$.
In addition, 
\begin{itemize}
    \item[$\cdot$]
    $N(\Pc_G) = 3  \cdot 2^{\frac{n-1}{2}}-2$
if and only if $G=K_{(n-1)/2,(n+1)/2}$.
    \item[$\cdot$]
$N(\Pc_G) =6^{\frac{n-1}{2}}$ if and only if $G$ is the $1$-sum of $(n-1)/2$ triangles.
\end{itemize}
    \item[(2)]
    If $n$ is even, then we have
$2^{\frac{n}{2}+1}-2 \le N(\Pc_G) \le 14 \cdot 6^{\frac{n}{2}-2}$.
In addition, 
\begin{itemize}
    \item[$\cdot$]
    $N(\Pc_G) = 2^{\frac{n}{2}+1}-2 $
if and only if $G=K_{n/2,n/2}$.
    \item[$\cdot$]
$N(\Pc_G) =14 \cdot 6^{\frac{n}{2}-2}$ if and only if $G$ is the $1$-sum  of $K_4$ with $n/2 -2$ triangles.
\end{itemize}

\end{itemize}
\end{Conjecture}

Let $G=(V,E)$ be a graph on the vertex set $V=[n-1]$.
Then the \textit{suspension} $\Ghat$ of $G$ is the graph on the vertex set 
$[n]$ and the edge set $E \cup \{\{i,n\} : i \in [n-1] \}$.
In the present paper, we show that Conjecture~\ref{yoso} is true for any suspension graph.

\begin{Theorem} \label{mainT}
    
Let $G$ be a graph on the vertex set $[n-1]$ with $n \ge 2$.
   Then 
   $$
   N(\Pc_{\Ghat}) \ge 2^{n-1}
   $$
    and equality holds if and only if $G$ is an empty graph (i.e., a graph having no edges), and hence $\Ghat$ is a star graph $K_{1,n-1}$.
    Moreover, 
   $$
   N(\Pc_{\Ghat}) \le 
   \begin{cases}
       6^{\frac{n-1}{2}} & \mbox{if } n \mbox{ is odd,}\\
       14 \cdot 6^{\frac{n}{2}-2} & \mbox{if } n \mbox{ is even}
                 \end{cases}
   $$
    and equality holds if and only if one of the following holds{\rm :}
    \begin{itemize}
        \item[{\rm (a)}]
        $n$ is odd, and
        $G$ is a disjoint union of $(n-1)/2$ edges,
        and hence $\Ghat$ is a $1$-sum of $(n-1)/2$ triangles.
        \item[{\rm (b)}]
        $n$ is even, and
         $G$ is a disjoint union of $n/2-2$ edges with a triangle,
         and hence $\Ghat$ is a $1$-sum of $K_4$ with $n/2-2$ triangles. 
    \end{itemize}
          
\end{Theorem}

In addition, we extend Theorem~\ref{mainT} to the join of two graphs.
Let $G_1 =(V, E)$ and $G_2=(V', E')$ be (not necessarily connected) graphs with $V \cap V' = \emptyset$.
Then the {\it join} $G_1 + G_2$ of $G_1$ and $G_2$ is the graph 
on the vertex set $V \cup V'$ and the edge set $E \cup E' \cup \{\{i,j\} : i \in V , j \in V'\}$.
For example, $K_\ell + K_m = K_{\ell+m}$ and the join of two empty graphs is a complete bipartite graph.
Note that $K_1 + G$ is the suspension of $G$.
By the following theorem, Conjecture \ref{yoso} holds for any 
connected graph whose complement is not connected.

\begin{Theorem} \label{JOIN theorem}
    Let $G_1 =(V, E)$ and $G_2=(V', E')$ be graphs with $V \cap V' = \emptyset$
    and let $n=|V| + |V'|$.
Then 
\[
3  \cdot 2^{\frac{n-1}{2}}-2 \le N(\Pc_{G_1 + G_2}) \le       6^{\frac{n-1}{2}}
\]
if $n$ is odd, and
\[
2^{\frac{n}{2}+1}-2 \le N(\Pc_{G_1 + G_2}) \le       14 \cdot 6^{\frac{n}{2}-2}
\]
if $n$ is even.
\end{Theorem}

The present paper is organized as follows.
In Section \ref{sec:two}, after reviewing the characterizations of the facets of symmetric edge polytopes,
we confirm that, in order to study Conjecture \ref{yoso},
it is enough to consider 
2-connected nonbipartite graphs.
Next, in Section \ref{sec:three}, using a characterization of the facets of
symmetric edge polytopes of suspension graphs,
we give a proof of Theorem~\ref{mainT}.
Finally, in Section \ref{sec:four}, we extend Theorem~\ref{mainT} to join graphs
by giving a proof of Theorem~\ref{JOIN theorem}.
From the results in the present paper, in order to study Conjecture \ref{yoso},
it is enough to discuss $2$-connected non-bipartite graphs whose complement is connected.

\section{Basics on the facets of symmetric edge polytopes}
\label{sec:two}

In the present section, 
we will give some basic results on the facets of symmetric edge polytopes.
First, we review the 
characterizations of facets of symmetric edge polytopes.
Let $G$ be a graph on the vertex set $V=[n]$.
A {\it spanning subgraph} of $G$ is a subgraph of $G$ which contains every vertex of $G$.
Since $\Pc_G$ is reflexive, it is known that
the supporting hyperplane of each facet of $\Pc_G$ is of the form 
$H=\{ \xb \in \RR^n : \langle \ab , \xb \rangle = 1\}$
for some vertex $\ab  \in \ZZ^n$ of $\Pc_G^\vee$.
By regarding $\ab = (a_1,\dots, a_n) \in \ZZ^n$ as the map
$f : V \rightarrow \ZZ$, $i \mapsto a_i$, we have the following.

\begin{Proposition}[{\cite[Theorem 3.1]{HJM}}] \label{facet defining}
Let $G = (V, E)$ be a connected graph. 
Then $f : V \rightarrow \ZZ$ defines a facet of $\Pc_G$ if and only if both of the following hold.
\begin{itemize}
\item[(i)]
For every edge $e = \{i,j\}\in E$, we have $|f(i) - f(j)| \le 1$.
\item[(ii)] 
The subset of edges $E_f := \{e = \{i,j\} \in  E : |f(i) - f(j)| = 1\}$
 forms a spanning connected subgraph of $G$.
\end{itemize}
 
\end{Proposition}

There exists a characterization for the subgraphs appearing in Proposition~\ref{facet defining}.

\begin{Definition}
If $f : V \rightarrow \ZZ$ defines a facet of $\Pc_G$, then 
the graph $G_f := (V, E_f)$ in Proposition~\ref{facet defining} is called the {\it facet subgraph} of $G$ associated with $f$.   
Let ${\rm FS}(G)$ denote the set of all facet subgraphs of $G$.
Given a facet subgraph $H \in  {\rm FS}(G) $, let $\mu(H)$ denote the number of facets of $\Pc_G$
whose facet subgraph is $H$.
\end{Definition}

Note that, if $G$ is bipartite, then ${\rm FS}(G) = \{G\}$.
The following fact is often used in the study of $N(\Pc_G)$.

\begin{Proposition}\label{fs sum}
    Let $G$ be a connected graph. Then 
\[
N(\Pc_G) = \sum_{H \in {\rm FS}(G)} \mu(H).
\]
\end{Proposition}

On the other hand, a characterization of facet subgraphs of $G$ is known.

\begin{Proposition}[{\cite[Theorem 3 (2)]{CDK}}] \label{mcsb}
Let $G$ be a connected graph.
A subgraph $H$ of $G$ is a facet subgraph of $G$
if and only if it is a maximal connected spanning bipartite subgraph of $G$.
\end{Proposition}

Let $G = (V, E)$ be a connected graph and let $H = (V, E_f)$ be a facet subgraph of $G$ associated with $f : V \rightarrow \ZZ $.
From Proposition \ref{facet defining}, we have $f(i) = f(j)$
for all $e = \{i, j\} \in E \setminus E_f$.
If the graph obtained by contracting all edges in $E \setminus E_f$ of $G$ and simplifying it (i.e., removing loops and multiple edges) is denoted by $G^*$, then we have $\mu(H) = N(\Pc_{G^*})$.
Given a subset $E' \subset E$,
let $G / E'$ denote the graph obtained by contracting all edges in $E'$ and simplifying it.

\begin{Example}
Let $G = ([5], E)$ be the following connected graph:

\begin{figure}[h]
\centering
\begin{tikzpicture}[scale=0.6]
\foreach \i in {1,...,5} {
    \coordinate (v\i) at ({90+72*(\i-1)}:2);
}

\draw[thick] (v1) -- (v2) node[midway, left] {$e_1$};
\draw[thick] (v2) -- (v3) node[midway, left] {$e_2$};
\draw[thick] (v3) -- (v4) node[midway, below] {$e_3$};
\draw[thick] (v4) -- (v5) node[midway, right] {$e_4$};
\draw[thick] (v5) -- (v1) node[midway, right] {$e_5$};
\draw[thick] (v2) -- (v4) node[midway, above] {$e_6$};
\draw[thick] (v3) -- (v5) node[midway, right] {$e_7$};

\foreach \i in {1,...,5} {
    \filldraw[black] (v\i) circle (3pt);
}
\end{tikzpicture}
\label{5cyle}
\end{figure}

\noindent
From Proposition \ref{mcsb}, ${\rm FS}(G) = \{H_1,H_2,\dots,H_7\}$, where the edge sets of each $H_k$ are given respectively by  
$$
E \setminus \{e_3\}, \;\; E \setminus \{e_2, e_7\}, \;\; E \setminus \{e_4, e_6\}, 
$$
$$
E \setminus \{e_1, e_2, e_4\}, \;\; E \setminus \{e_1, e_6, e_7\}, \;\; E \setminus \{e_2, e_4, e_5\}, \;\; E \setminus \{e_5, e_6, e_7\}.
$$    
Then $G / \{e_3\}$ is isomorphic to $K_{2,2}$,
$G / \{e_2, e_7\}$ and $G / \{e_4, e_6\}$ are isomorphic to $K_{1,2}$,
$G / \{e_1, e_2, e_4\}$, $G / \{e_1, e_6, e_7\}$, $G / \{e_2, e_4, e_5\}$ and $G / \{e_5, e_6, e_7\}$ are isomorphic to $K_{1,1}$, respectively.
By Proposition \ref{fs sum} and
(\ref{2bu graph}) in Introduction,
\[
N(\Pc_G) = \sum_{k=1}^7 \mu(H_k) = N(\Pc_{K_{2,2}}) + 2 N(\Pc_{K_{1,2}}) + 4 N(\Pc_{K_{1,1}}) = 22
\] 
holds.
\end{Example}

The following upper bound for bipartite graphs is known.

\begin{Proposition}[{\cite[Corollary 33]{DDM}}]\label{tree facet}
    Let $G$ be a connected bipartite graph with $n$ vertices.
Then $N(\Pc_G) \leq 2^{n-1}$, and the equality holds if $G$ is a tree.
\end{Proposition}

Note that $2^{n-1} < 14 \cdot 6^{\frac{n}{2}-2}$ for any $n \ge 2$.
Since we cannot find it in literature, 
we confirm that the lower bound in Conjecture~\ref{yoso} is true for
bipartite graphs
by using the following proposition.

\begin{Proposition} \label{hikinukiB}
Let $G=(V,E)$ be a connected bipartite graph.
Suppose that the bipartite graph $G-e$ on the vertex set $V$ obtained from $G$ by deleting an edge $e$ of $G$ 
is connected.
Then we have
$N(\Pc_G) \le N(\Pc_{G - e})$.
\end{Proposition}

\begin{proof}
Since both $G$ and $G -e$ are bipartite,
we have $N(\Pc_G) = \mu (G)$ and $N(\Pc_{G-e}) = \mu (G-e)$
from Propositions \ref{fs sum} and \ref{mcsb}.
From Proposition~\ref{facet defining}, 
it follows that $N(\Pc_G) = \mu(G) \le \mu(G-e) = N(\Pc_{G - e})$.
\end{proof}

We now show that the lower bound in Conjecture~\ref{yoso} is true for bipartite graphs.

\begin{Proposition} \label{bipartite conj}
Let $G$ be a connected bipartite graph with $n$ vertices. 
\begin{itemize}
\item[(a)]
If $n$ is odd, then we have $N(\Pc_G) \geq  3  \cdot 2^{\frac{n-1}{2}}-2$, and the equality holds if and only if $G=K_{(n-1)/2,(n+1)/2}$.
    \item[(b)]
If $n$ is even, then we have $N(\Pc_G) \geq 2^{\frac{n}{2}+1}-2$, and the equality holds if and only if $G=K_{n/2,n/2}$.
\end{itemize}
\end{Proposition}

\begin{proof}
Let $G$ be a connected bipartite graph on the vertex set $V = V_1 \sqcup V_2$, where $n_1 = |V_1|$ and $n_2 = |V_2|$
with $n_1 \le n_2$.
Equation (\ref{2bu graph}) in Introduction
says that 
$N(\Pc_{K_{n_1,n_2}}) = 2^{n_1}+2^{n_2}-2$.
Using Proposition \ref{hikinukiB} repeatedly
from $K_{n_1,n_2}$ to $G$, we have 
\[
 N(\Pc_G) 
\ge
N(\Pc_{K_{n_1,n_2}}) 
=
2^{n_1}+2^{n_2}-2. 
\]
On the other hand, 
$$
2^{n_1}+2^{n_2}
\ge
\left\{
\begin{array}{cc}
2^{\frac{n}{2}} + 2^{\frac{n}{2}} =2^{\frac{n}{2}+1} & \mbox{if } n \mbox{ is even,}\\
2^{\frac{n-1}{2}} + 2^{\frac{n+1}{2}} = 3 \cdot 2^{\frac{n-1}{2}} & \mbox{if } n \mbox{ is odd.}
\end{array}
\right.
$$
Equality holds when $n/2=n_1=n_2$ if $n$ is even,
and $n_1 = (n-1)/2$ and $n_2 = (n+1)/2$ if $n$ is odd.
\end{proof}

We close the present section by proving that a connected graph $G$ satisfies 
the condition of Conjecture~\ref{yoso} if 
each ``block" of $G$ satisfies the condition of Conjecture~\ref{yoso}.
Blocks of a graph are defined as follows.

\begin{Definition}
Let $G$ be a connected graph.
A vertex $v$ of $G$ is called a {\it cut vertex} if the graph obtained by
the removal of $v$ from $G$ is disconnected.
A \textit{block} of $G$ is a maximal connected subgraph of $G$ without cut vertices.
\end{Definition}

In particular, any connected graph is the $1$-sum of its blocks.

\begin{Proposition}[{\cite[Proposition 9]{BrBr}}] \label{0sum}
Let $G$ be the $1$-sum of connected graphs $G_1$ and $G_2$.
Then we have $N(\Pc_G) =  N(\Pc_{G_1})N(\Pc_{G_2})$.    
\end{Proposition}

From this proposition, we have the following.

\begin{Proposition} \label{block prop}
Let $G$ be the $1$-sum of connected graphs $G_1$ and $G_2$.
If $G_1$ and $G_2$ satisfy the condition of Conjecture~\ref{yoso},
then so does $G$.
  
\end{Proposition}
    
\begin{proof}
Let $n_i \ge 2$ be the number of vertices of $G_i$ for $i=1,2$.
Then $G$ has $n=n_1 +n_2-1$ vertices.

\medskip

\noindent
{\bf Case 1} (both $n_1$ and $n_2$ are odd){\bf .}
Then $n$ is odd. From Proposition~\ref{0sum},
\[
 N(\Pc_G) = N(\Pc_{G_1})N(\Pc_{G_2})
 \le 
 6^{\frac{n_1-1}{2}} \cdot  6^{\frac{n_2-1}{2}}
 =6^{\frac{n-1}{2}}
\]
and
\begin{eqnarray*}
 N(\Pc_G) = N(\Pc_{G_1})N(\Pc_{G_2})
 &\ge& 
\left(3  \cdot 2^{\frac{n_1-1}{2}}-2\right)
\left(3  \cdot 2^{\frac{n_2-1}{2}}-2\right)\\
 &=& 3  \cdot 2^{\frac{n-1}{2}} -2
+ 6 \left(2^{\frac{n_1-1}{2}}-1\right)\left(2^{\frac{n_2-1}{2}}-1\right)\\
& > & 3  \cdot 2^{\frac{n-1}{2}} -2.
\end{eqnarray*}

\medskip

 \noindent
{\bf Case 2} (both $n_1$ and $n_2$ are even){\bf .}
Then $n$ is odd.
From Proposition~\ref{0sum},
\[
 N(\Pc_G) = N(\Pc_{G_1})N(\Pc_{G_2})
 \le 
 14 \cdot 6^{\frac{n_1}{2}-2} \cdot  14 \cdot 6^{\frac{n_2}{2}-2}
 =\frac{49}{54} \cdot 6^{\frac{n-1}{2}}
 <  6^{\frac{n-1}{2}}
\]
and
\begin{eqnarray*}
 N(\Pc_G) = N(\Pc_{G_1})N(\Pc_{G_2})
 &\ge& 
\left(2^{\frac{n_1}{2}+1}-2 \right)
\left(2^{\frac{n_2}{2}+1}-2 \right)\\
 &=& 3  \cdot 2^{\frac{n-1}{2}} -2+
 2\left(2^{\frac{n_1}{2}}-2 \right)
\left(2^{\frac{n_2}{2}}-2 \right)+
 2^{\frac{n-1}{2}}-2
 \\
& \ge & 3  \cdot 2^{\frac{n-1}{2}} -2.
\end{eqnarray*}
(In the last inequality, equality holds if and only if $n_1 = n_2 =2$ and hence $G=K_{1,2}$.)

\medskip

\noindent
{\bf Case 3} ($n_1$ is odd and $n_2$ is even){\bf .}
Then $n$ is even.
From Proposition~\ref{0sum},
\[
 N(\Pc_G) = N(\Pc_{G_1})N(\Pc_{G_2})
 \le 
 6^{\frac{n_1-1}{2}} \cdot  14 \cdot 6^{\frac{n_2}{2}-2}
 =14 \cdot 6^{\frac{n}{2}-2}
\]
and
\begin{eqnarray*}
 N(\Pc_G) = N(\Pc_{G_1})N(\Pc_{G_2})
 &\ge& 
\left(3  \cdot 2^{\frac{n_1-1}{2}}-2 \right)
\left(2^{\frac{n_2}{2}+1}-2 \right)\\
 &=& 2^{\frac{n}{2}+1}-2+
 2 \left(  2^{\frac{n_1-1}{2}} -1\right)
 \left( 2^{\frac{n_2}{2}+1}-3\right)
 \\
& > & 2^{\frac{n}{2}+1}-2,
\end{eqnarray*}
as desired.
\end{proof}

As explained in Introduction, 
it is known that $N(\Pc_{K_{\ell,m}}) = 2^\ell + 2^m-2$ and
$N(\Pc_{K_{\ell_1,\dots,\ell_s}}) = 2^{\sum_{i=1}^s \ell_i} -  \sum_{i=1}^s ( 2^{\ell_i} -2) -2$ if $s \ge 3$.
Thus Conjecture~\ref{yoso} is true for complete multipartite graphs.
Since every 2-connected graph with $n \le 4$ vertices is complete multipartite,
we have the following from Proposition~\ref{block prop}.

\begin{Proposition} \label{small}
    Conjecture \ref{yoso} is true for $n=3,4$.
\end{Proposition}

\section{Facets of symmetric edge polytopes of suspension graphs}
\label{sec:three}

In the present section, 
using a characterization of the facets of
symmetric edge polytopes of suspension graphs,
we give a proof of Theorem~\ref{mainT}.

\begin{Definition}
Let $G$ be a graph on the vertex set $V$.
Given a vertex $v$ of $G$, let $N_G(v)$ denote the set of all vertices
that are adjacent to $v$ in $G$.
Let $N_G[v] := N_G(v) \cup \{v\}$.
A subset $S \subset V$ is called a {\it dominating set} of $G$ if 
$\bigcup_{v \in S} N_G[v] = V$.
\end{Definition}

Note that if $S \subset V$ is a dominating set of $G$, then any $S' \subset V$ with $S \subset S'$ is a
dominating set of $G$.
Facet subgraphs of a suspension graph is characterized by dominating sets.

\begin{Lemma} \label{domi}
Let $G$ be a graph on the vertex set $[n-1]$,
and let $H$ be a maximal spanning bipartite subgraph of $\Ghat$ on the vertex set $[n]=V_1 \sqcup V_2$, where $n \in V_1$.
Then $H$ is a facet subgraph of $\Ghat$ if and only if $V_2$ is a dominating set of $G$.
\end{Lemma}

\begin{proof}
    Since $H$ is a maximal spanning bipartite subgraph of $\Ghat$, $H$ is a facet subgraph of $\Ghat$ if and only if 
    $H$ is connected.
    Since $n \in V_1$ is adjacent to any vertex in $V_2$, $H$ is connected if and only if $V_2$ is a dominating set of $G$.
\end{proof}

\begin{Definition}
Let $G$ be a graph on the vertex set $V$.
Then let $c(G)$ denote the number of connected components of $G$.    
Given a subset $S \subset V$, let $G[S]$ denote the induced subgraph 
of $G$ on the vertex set $S$.
\end{Definition}

\begin{Lemma} \label{countlemma}
    Let $G$ be a graph on the vertex set $[n-1]$.
    Suppose that $H$ is a facet subgraph of $\Ghat$ on the vertex set $[n] = V_1 \sqcup V_2$, where $n \in V_1$.
    Then we have $\mu(H) = 2^{c(G[V_2])}$.
\end{Lemma}

\begin{proof}
Suppose that $H$ is the facet subgraph for a facet defined by $f : [n] \rightarrow \ZZ$.
We may assume that $f(n)=0$.
For each $i \in V_1$, since
$\{i,n\}$ is an edge of $\widehat{G}$ and not an edge of $H$, 
we have $f(i)=0 $ from Proposition~\ref{facet defining}.
Since $H$ is a facet subgraph of $f$, it follows that
$|f(j)|=1$ for each $j \in V_2$.
If 
$j_1, j_2 \in V_2$ belong to the same connected component of $G[V_2]$, then $f(j_1) = f(j_2)$.
If $j_1, j_2 \in V_2$ do not belong to the same connected component of $G[V_2]$, then $f(j_1)$ and $f(j_2)$ are independent.
Thus one can choose $1$ or $-1$ for the value of $f$ 
for each connected component of $G[V_2]$.
 \end{proof}

 We give an example of $\mu(H)$ from Lemma \ref{countlemma} in Figure \ref{ExampleOfLemma}.
The graph obtained by contracting all edges in $\{\{i, j\} : i, j \in V_1\} \cup \{\{i, j\} : i, j \in V_2\}$ of $\Ghat$ (and simplifying it) is isomorphic to $K_{1,3}$.
Hence, $\mu(H) = N(\Pc_{K_{1, 3}}) = 2^3$ holds.

\begin{figure}
\begin{center}
\scalebox{0.9}{
\begin{tikzpicture}
\coordinate (n0) at (0,1.3) node at (n0) {$n$};
\coordinate(n1) at (-0.5,1) node at (n1) {$V_1$};
\coordinate (n2) at (-0.5,-1) node at (n2)  {$V_2$};
\coordinate (n3) at (2,-1.5) node at (n3) {$\Ghat$};

\coordinate(v1) at (0,1);
\coordinate (v2) at (1,1);
\coordinate (v3) at (2,1);
\coordinate (v4) at (3,1);
\coordinate (v5) at (4,1);

\coordinate (w1) at (0,-1);
\coordinate (w2) at (1,-1);
\coordinate (w3) at (2,-1);
\coordinate (w4) at (3,-1);
\coordinate (w5) at (4,-1);


\draw 
(v1)--(w1)
(v1)--(w2)
(v1)--(w3)
(v1)--(w4)
(v1)--(w5)
(v2)--(w1)
(v3)--(w2)
(v3)--(w3)
(v4)--(w1)
(v4)--(w5)
(v5)--(w4)
(v1) to [out=30,in=150] (v2)
(v1) to [out=35,in=145] (v3)
(v1) to [out=40,in=140] (v4)
(v1) to [out=45,in=135] (v5)
(v3) to [out=40,in=140] (v5)
(w2) to [out=-30,in=-150] (w3)
(w4) to [out=-30,in=-150] (w5)
;
\fill
(v1) circle  (2pt)
(v2) circle  (2pt)
(v3) circle  (2pt)
(v4) circle  (2pt)
(v5) circle  (2pt)
(w1) circle  (2pt)
(w2) circle  (2pt)
(w3) circle  (2pt)
(w4) circle  (2pt)
(w5) circle  (2pt)
;

\draw [very thick, -latex] (5,0) -- (6,0);


\coordinate(v11) at (7,1) node at (v11) [above] {$0$};
\coordinate (v12) at (8,1)node at (v12) [above] {$0$};
\coordinate (v13) at (9,1) node at (v13) [above] {$0$};
\coordinate (v14) at (10,1) node at (v14) [above] {$0$};
\coordinate (v15) at (11,1) node at (v15) [above] {$0$};

\coordinate (w11) at (7,-1) node at (w11) [below] {$1$};
\coordinate (w12) at (8,-1)node at (w12) [below] {$-1$};
\coordinate (w13) at (9,-1)node at (w13) [below] {$-1$};
\coordinate (w14) at (10,-1)node at (w14) [below] {$-1$};
\coordinate (w15) at (11,-1)node at (w15) [below] {$-1$};

\draw 
(v11)--(w11)
(v11)--(w12)
(v11)--(w13)
(v11)--(w14)
(v11)--(w15)
(v12)--(w11)
(v13)--(w12)
(v13)--(w13)
(v14)--(w11)
(v14)--(w15)
(v15)--(w14)
(v11) to [out=30,in=150] (v12)
(v11) to [out=35,in=145] (v13)
(v11) to [out=40,in=140] (v14)
(v11) to [out=45,in=135] (v15)
(v13) to [out=40,in=140] (v15)
(w12) to [out=-30,in=-150] (w13)
(w14) to [out=-30,in=-150] (w15)
;
\fill
(v11) circle  (2pt)
(v12) circle  (2pt)
(v13) circle  (2pt)
(v14) circle  (2pt)
(v15) circle  (2pt)
(w11) circle  (2pt)
(w12) circle  (2pt)
(w13) circle  (2pt)
(w14) circle  (2pt)
(w15) circle  (2pt)
;

\draw [very thick, -latex] (12,0) -- (13,0);

\coordinate (v21) at (15,1) node at (v21) [above] {$0$};

\coordinate (w21) at (14,-1) node at (w21) [below] {$1$};
\coordinate (w22) at (15,-1) node at (w22) [below] {$-1$};
\coordinate (w23) at (16,-1) node at (w23) [below] {$-1$};

\draw 
(v21)--(w21)
(v21)--(w22)
(v21)--(w23)
;

\fill
(v21) circle (2pt)
(w21) circle (2pt)
(w22) circle (2pt)
(w23) circle (2pt)
;

\end{tikzpicture}
}
\end{center}
\caption{Example of $\mu(H)$ from Lemma \ref{countlemma}.}\label{ExampleOfLemma}
\end{figure}

\begin{Definition}
Given a vertex $v$ of a graph $G=(V,E)$, we define the following three graphs:
\begin{itemize}
    \item
    Let $G-v$ denote the induced subgraph $G[V \setminus \{v\}]$ of $G$;
    \item 
    If $N_G[v] \ne V$, then let $G - N_G[v]$ denote the induced subgraph $G[V \setminus N_G[v]]$ of $G$;
    \item 
    Let $G/v$ denote the graph obtained from $G$ by removal of $v$
and insertion of all edges $\{i,j\}$ such that $i ,j \in N_G(v)$.

\end{itemize}
\end{Definition}

For example, if $G$ is a graph
and $v$ is a vertex of $G$ as in Figure~\ref{Gobr}, then $G-v$, $G-N_G(v)$ and $G/v$ are graphs in Figure~\ref{Gobr}.

\begin{figure}
    \centering
    \begin{tikzpicture}
\coordinate(v1) at (0,1) node at (v1) [above]{$v$};
\coordinate(v2) at (-0.95,0.30);
\coordinate(v3) at (-0.58,-0.80);
\coordinate(v4) at (0.58,-0.80);
\coordinate (v5) at (0.95,0.30);
\coordinate (v6) at (0,0);
\coordinate (n1) at (0,-1) node at (n1) [below] {$G$};

\draw
(v1)--(v2)--(v3)--(v4)--(v5)--(v1)
(v1)--(v6)
(v2)--(v6)
(v3)--(v6)
(v4)--(v6)
(v5)--(v6)
;

\fill
(v1) circle (1.8pt)
(v2) circle (1.8pt)
(v3) circle (1.8pt)
(v4) circle (1.8pt)
(v5) circle (1.8pt)
(v6) circle (1.8pt)
;

\coordinate(v2) at (-0.95+3,0.30);
\coordinate(v3) at (-0.58+3,-0.80);
\coordinate(v4) at (0.58+3,-0.80);
\coordinate (v5) at (0.95+3,0.30);
\coordinate (v6) at (0+3,0);
\coordinate (n1) at (0+3,-1) node at (n1) [below] {$G-v$};

\draw
(v2)--(v3)--(v4)--(v5)
(v2)--(v6)
(v3)--(v6)
(v4)--(v6)
(v5)--(v6)
;

\fill
(v2) circle (1.8pt)
(v3) circle (1.8pt)
(v4) circle (1.8pt)
(v5) circle (1.8pt)
(v6) circle (1.8pt)
;

\coordinate(v3) at (-0.58+6,-0.80);
\coordinate(v4) at (0.58+6,-0.80);
\coordinate (n1) at (6,-1) node at (n1) [below] {$G-N_G[v]$};

\draw
(v3)--(v4)
;

\fill
(v3) circle (1.8pt)
(v4) circle (1.8pt)
;

\coordinate(v2) at (-0.95+8.75,0.30);
\coordinate(v3) at (-0.58+8.75,-0.80);
\coordinate(v4) at (0.58+8.75,-0.80);
\coordinate (v5) at (0.95+8.75,0.30);
\coordinate (v6) at (0+8.75,0);
\coordinate (n1) at (0+8.75,-1) node at (n1) [below] {$G/v$};

\draw
(v2)--(v3)--(v4)--(v5)--(v6)--(v2)
(v2)--(v6)
(v3)--(v6)
(v4)--(v6)
(v5)--(v6)
(v2)--(v5)
;

\fill
(v1) circle (1.8pt)
(v2) circle (1.8pt)
(v3) circle (1.8pt)
(v4) circle (1.8pt)
(v5) circle (1.8pt)
(v6) circle (1.8pt)
;

\end{tikzpicture}
    \caption{Three graphs obtained from $G$ and $v$.}
    \label{Gobr}
\end{figure}

\begin{Proposition}
\label{insp}
    Let $G$ be a graph on the vertex set $[n-1]$ with $n \ge 3$.
    Given a vertex $v$ of $G$,
    we have 
    \begin{equation} \label{keyF}
    N(\Pc_{\widehat{G-v}})
    + N(\Pc_{\widehat{G/v}}) \ \ 
    \le \ \ N(\Pc_\Ghat) \ \  \le  \ \ 
    N(\Pc_{\widehat{G-v}})+ 2N(\Pc_{\widehat{G-N_G[v]}})+ N(\Pc_{\widehat{G/v}})
    \end{equation}
    if $N_G[v] \ne [n-1]$, and
        \begin{equation} \label{keyF2}
   N(\Pc_\Ghat) \ \  =  \ \ 
    N(\Pc_{\widehat{G-v}})+  N(\Pc_{\widehat{G/v}}) +2
    \end{equation}
    if $N_G[v] = [n-1]$.
    \end{Proposition} 

\begin{proof}
From Proposition \ref{fs sum}, we have
$$
\begin{array}{cc} 
N(\Pc_\Ghat)=\sum_{H \in {\rm FS}(\Ghat)} \mu(H),&
N(\Pc_{\widehat{G-v}})=\sum_{H \in {\rm FS}({\widehat{G-v}})} \mu(H),\\
N(\Pc_{\widehat{G-N_G[v]}})=\sum_{H \in {\rm FS}({\widehat{G-N_G[v]}})} \mu(H),&
N(\Pc_{\widehat{G/v}})=\sum_{H \in {\rm FS}({\widehat{G/v}})} \mu(H).
\end{array}
$$
In order to compare these values, we define partitions
\begin{eqnarray*}
{\rm FS} (\widehat{G}) &=& {\rm FS}_0 (\widehat{G}) \sqcup {\rm FS}_1 (\widehat{G}) \sqcup {\rm FS}_2 (\widehat{G}),\\
{\rm FS} (\widehat{G-v}) &=& {\rm FS}_1 (\widehat{G-v}) \sqcup {\rm FS}_2 (\widehat{G-v}),\\
{\rm FS} (\widehat{G/v}) &=& {\rm FS}_1 (\widehat{G/v}) \sqcup {\rm FS}_2 (\widehat{G/v}),
\end{eqnarray*}
where 
\begin{eqnarray*}
{\rm FS}_0 (\widehat{G}) &:=&
\left\{ H \in {\rm FS}(\Ghat) : 
\mbox{the bipartition of } H \mbox{ is } V_1 \sqcup V_2, \mbox{ where }
n , v \in V_1
\right\},\\
{\rm FS}_1 (\widehat{G}) &:=& \left\{ H \in {\rm FS}(\Ghat)  :
\begin{array}{c} \mbox{the bipartition of } H \mbox{ is } V_1 \sqcup V_2, \mbox{ where }
n \in V_1, v \in V_2 \\\ 
\mbox{ and }  N_G(v) \subset V_1
\end{array} 
\right\},\\
{\rm FS}_2 (\widehat{G}) &:=& \left\{ H \in {\rm FS}(\Ghat)  :
\begin{array}{c} \mbox{the bipartition of } H \mbox{ is } V_1 \sqcup V_2, \mbox{ where }
n\in V_1, v \in V_2\\
 \mbox{ and } N_G(v) \cap V_2 \ne \emptyset
\end{array} 
\right\},\\
{\rm FS}_1 (\widehat{G-v}) &:=& \left\{ H \in {\rm FS}(\widehat{G-v})  :
\begin{array}{c} \mbox{the bipartition of } H \mbox{ is } V_1 \sqcup V_2, \mbox{ where }
n \in V_1, \\
\mbox{and } N_G(v) \subset V_1 \end{array}\right\},\\
{\rm FS}_2 (\widehat{G-v}) &:=& \left\{ H \in {\rm FS}(\widehat{G-v})  :
\begin{array}{c} \mbox{the bipartition of } H \mbox{ is } V_1 \sqcup V_2, \mbox{ where }
n \in V_1, \\
\mbox{and } N_G(v) \cap V_2 \ne \emptyset \end{array}\right\},\\
{\rm FS}_1 (\widehat{G/v}) &:=& \left\{ H \in {\rm FS}(\widehat{G/v})  :
\begin{array}{c} \mbox{the bipartition of } H \mbox{ is } V_1 \sqcup V_2, \mbox{ where }
n \in V_1, \\
\mbox{and } N_G(v) \subset V_1 \end{array}\right\},\\
{\rm FS}_2 (\widehat{G/v}) &:=& \left\{ H \in {\rm FS}(\widehat{G/v})  :
\begin{array}{c} \mbox{the bipartition of } H \mbox{ is } V_1 \sqcup V_2, \mbox{ where }
n \in V_1, \\
\mbox{and } N_G(v) \cap V_2 \ne \emptyset \end{array}\right\}.
\end{eqnarray*}

First, we will show the following equalities and inequalities:
\begin{eqnarray}
    \sum_{H \in {\rm FS}_0(\Ghat) } \mu (H) &=& \sum_{H \in {\rm FS}_2(\widehat{G-v})} \mu(H),\label{c1}\\ 
    \sum_{H \in {\rm FS}_2(\Ghat) } \mu (H) &=& \sum_{H \in {\rm FS}_2(\widehat{G/v})} \mu(H),\label{c2}\\ 
    \sum_{H \in {\rm FS}_1(\Ghat) } \mu (H) & & 
    \left\{
    \begin{array}{ccc}
     \le &  \displaystyle 2 \sum_{H \in  {\rm FS} (\widehat{G-N_G[v]})} \mu(H)
     &  \mbox{if } N_G[v] \ne [n-1],\\
     \\
     = &2    &  \mbox{if } N_G[v] = [n-1],
    \end{array}
    \right.\label{c3}\\ 
    \sum_{H \in {\rm FS}_1(\widehat{G-v}) } \mu (H) &\le& \frac{1}{2} \sum_{H \in {\rm FS}_1(\Ghat)} \mu(H),\label{c4}\\ 
    \sum_{H \in {\rm FS}_1(\widehat{G/v}) } \mu (H) &\le& \frac{1}{2}\sum_{H \in {\rm FS}_1(\Ghat)} \mu(H). \label{c5}
\end{eqnarray}

\noindent
{\bf Proof of (\ref{c1}).}
It is enough to show that
$\varphi: {\rm FS}_0(\Ghat) \rightarrow {\rm FS}_2(\widehat{G-v})$,
$H \mapsto H-v$ is a bijection
such that $\mu(H) = \mu(\varphi(H))$.
%

Let $H \in {\rm FS}_0(\Ghat)$ and let $[n] = V_1 \sqcup V_2$
be the bipartition of the vertex set of $H$ with $n,v \in V_1$.
From Lemma \ref{domi}, $V_2$ is a dominating set of $G$.
In particular, $N_G(v) \cap V_2 \ne \emptyset$.
It is trivial that $V_2 \subset  \bigcup_{v' \in V_2} N_{G-v}[v']$.
Let $v_1 \in V_1 \setminus \{v\}$.
Since $V_2$ is a dominating set of $G$,
there exists $v_2 \in V_2$ such that
$\{v_1,v_2\}$ is an edge of $G$.
Since $v_1, v_2 \ne v$, $\{v_1,v_2\}$ is an edge of $G -v$.
Thus $v_1 \in \bigcup_{v' \in V_2} N_{G-v}[v']$
and hence $V_1 \setminus \{v\} \subset \bigcup_{v' \in V_2} N_{G-v}[v']$.
Therefore
$$[n]\setminus \{v\} =  (V_1 \setminus \{v\}) \sqcup V_2 = \bigcup_{v' \in V_2} N_{G-v}[v'],$$ i.e.,
$V_2$ is a dominating set of $G-v$.
From Lemma \ref{domi}, we have
$H - v \in {\rm FS}_2 (\widehat{G - v})$ 
since $H-v$ is a maximal spanning bipartite subgraph of $\widehat{G-v}$ on the vertex set $(V_1 \setminus \{v\}) \sqcup V_2$.
In addition, from $v \notin V_2$, we have $G[V_2] = (G-v) [V_2]$, 
and hence $c(G[V_2]) = c( (G-v)[V_2] )$.
From Lemma~\ref{countlemma}, 
$\mu(H) = \mu(H-v)=2^{c(G[V_2])}$.

Conversely, let $H_1 \in {\rm FS}_2 (\widehat{G-v})$
and let $[n] \setminus \{v\} = V_1' \sqcup V_2'$ be the bipartition of the vertex set of $H_1$
where $n \in V_1'$ and $N_G(v) \cap V_2' \ne \emptyset$.
Then $V_2'$ is a dominating set of $G-v$,
and hence $[n] \setminus \{v\} =  \bigcup_{v' \in V_2'} N_{G-v}[v'] \subset  \bigcup_{v' \in V_2'} N_{G}[v']$.
Moreover since $N_G(v) \cap V_2' \ne \emptyset$,
$v \in  \bigcup_{v' \in V_2'} N_{G}[v']$.
Thus $V_2'$ is a dominating set of $G$.
Hence the bipartite graph $H$ obtained 
from $H_1$ by adding the vertex $v$
and 
edges $\{v, v'\}$ ($v' \in N_G(v) \cap V_2'$)
is a facet subgraph of $\Ghat$
with $H \in {\rm FS}_0(\Ghat)$ and $\varphi(H)=H_1$.

\medskip

\noindent
{\bf Proof of (\ref{c2}).} 
It is enough to show that
$\varphi: {\rm FS}_2(\Ghat) \rightarrow {\rm FS}_2 (\widehat{G/v})$,
$H \mapsto \widetilde{H}$ defined below is a bijection
such that $\mu(H) = \mu(\varphi(H))$.
%

Let $H\in {\rm FS}_2(\Ghat)$
and let $[n]= V_1 \sqcup V_2$ be the bipartition of the vertex set of $H$ where $n\in V_1$, $v \in V_2$ and $N_G(v) \cap V_2 \neq \emptyset$.
Then $V_2$ is a dominating set of $G$.
It is trivial that $V_2 \setminus \{v\} \subset \bigcup_{v' \in V_2 \setminus \{v\}} N_{G/v}[v']$.
Let $v_1 \in V_1$.
Since $V_2$ is a dominating set of $G$, there exists an edge $\{v_1,v_2\}$ of $G$ for some $v_2 \in V_2$.
If $v_2 \neq v$, then $\{v_1,v_2\}$ is an edge of $G/v$, and hence $v_1 \in \bigcup_{v' \in V_2 \setminus \{v\}} N_{G/v}[v']$.
Suppose that $v_2 =v$.
Then $v_1 \in N_G(v)$.
Since $N_G(v) \cap V_2 \ne \emptyset$,
there exists $v_3 \in N_G(v) \cap V_2$.
Then $\{v_1,v_3\}$ is an edge of $G/v$, and hence $v_1 \in \bigcup_{v' \in V_2 \setminus \{v\}} N_{G/v}[v']$.
Thus $V_1 \subset \bigcup_{v' \in  V_2 \setminus \{v\}} N_{G/v}[v']$.
Therefore,
$$
[n] \setminus \{v\} = V_1 \sqcup (V_2 \setminus \{v\}) = \bigcup_{v' \in V_2 \setminus \{v\}} N_{G/v}[v']
$$
i.e., $V_2 \setminus \{v\}$ is a dominating set of $G/v$.
Hence 
 the maximal spanning bipartite subgraph $\widetilde{H}$ of $\widehat{G/v}$
on the vertex set $V_1 \sqcup (V_2 \setminus \{v\})$
is a facet subgraph of $\widehat{G/v}$.
Since $N_G(v) \cap V_2 \ne \emptyset$,
we have $\widetilde{H} \in {\rm FS}_2 (\widehat{G/v})$.
We now show that $c(G[V_2]) = c((G/v) [V_2 \setminus \{v\}])$.
Since $N_G(v) \cap V_2 \neq \emptyset$, $v$ is not an isolated vertex in $G[V_2]$.
Let $X_1,\dots, X_{t}$ with $t = c(G[V_2])$ be connected components of $G[V_2]$.
By definition, $N_G[v] \cap V_2 \subset X_i$ for some $i$.
Then $(G/v) [V_2 \setminus \{v\}]$ is obtained from $G[V_2]$ by removing $v$ from $X_i$
and inserting all edges $\{i,j\}$ such that $i ,j \in N_G(v) \cap V_2 \subset X_i$.
Thus $X_1,\dots, X_{i-1}, X_i \setminus \{v\} ,X_{i+1} ,\dots, X_{t}$ are connected components of $(G/v) [V_2 \setminus \{v\}]$,
and hence
$c(G[V_2]) = c((G/v) [V_2 \setminus \{v\}])$.
From Lemma~\ref{countlemma}, 
$\mu(H)= \mu(\widetilde{H}) = 2^{c(G[V_2])}$.

Conversely, let $H_2 \in {\rm FS}_2 (\widehat{G/v})$
and let $[n] \setminus \{v\}=V_1' \sqcup V_2'$ be the bipartition of the vertex set of $H_2$ where $n \in V_1'$ and $N_G(v) \cap V_2' \ne \emptyset$.
Then $V_2'$ is a dominating set of $G/v$.
Let $S:=V_2' \cup \{v\}$.
It is trivial that $S \subset \bigcup_{v' \in S} N_{G}[v']$.
Let $v_1' \in V_1'$.
Since $V_2'$ is a dominating set of $G/v$,
there exists $v_2' \in V_2'$ such that $\{v_1',v_2'\}$
is an edge of $G/v$.
If $\{v_1',v_2'\}$ is an edge of $G$, then
$v_1' \in N_{G}[v_2'] \subset \bigcup_{v' \in S} N_{G}[v']$.
Suppose that $\{v_1',v_2'\}$
is not an edge of $G$.
Then $v_1',v_2' \in N_G(v)\subset \bigcup_{v' \in S} N_{G}[v']$.
Hence we have $V_1' \subset  \bigcup_{v' \in S} N_{G}[v']$.
Thus
$$
[n] = V_1' \sqcup S  =  \bigcup_{v' \in S} N_{G}[v']
$$
i.e., $S$ is a dominating set of $G$.
Hence the maximal spanning bipartite subgraph $H$ of $\Ghat$ on the vertex set $V_1' \sqcup S$
is a facet subgraph of $\Ghat$
such that $\varphi(H) = H_2$.

\medskip

\noindent
{\bf Proof of (\ref{c3}).} Suppose that
$N_G[v] \ne [n-1]$.
Then it is enough to show that $\varphi: {\rm FS}_1(\Ghat) \rightarrow {\rm FS} (\widehat{G-N_G[v]})$,
 $H \mapsto H-N_G[v]$
 is an injection 
 such that $\mu(H) = 2 \mu(\varphi(H))$.


Let $H\in {\rm FS}_1(\Ghat)$
and let $[n]=V_1 \sqcup V_2$ be the bipartition of the vertex set of $H$ 
where $n \in V_1$, $v \in V_2$ and $N_G(v) \subset V_1$.
Then $V_2$ is a dominating set of $G$.
Let $S=V_2 \setminus \{v\}$.
It is trivial that $S \subset  \bigcup_{v' \in S} N_{G-N_G[v]}[v']$.
Let $v_1 \in V_1 \setminus N_G(v)$.
Since $V_2$ is a dominating set of $G$, $\{v_1,v_2\}$
is an edge of $G$ for some $v_2 \in V_2$.
Then $v \ne v_2$ from $v_1 \notin N_G(v)$.
Since $N_G(v) \subset V_1$, $v_2 \notin N_G[v]$.
Hence  $\{v_1,v_2\}$ is an edge of $G-N_G[v]$.
Thus $v_1 \in N_{G-N_G[v]}[v_2] \subset \bigcup_{v' \in S} N_{G-N_G[v]}[v']$.
It follows that $V_1 \setminus N_G(v) \subset  \bigcup_{v' \in S} N_{G-N_G[v]}[v']$.
Therefore 
$$
[n]\setminus N_G[v] =  (V_1 \setminus N_G(v)) \sqcup S =  \bigcup_{v' \in S} N_{G-N_G[v]}[v']
$$
i.e., $S$ is a dominating set of $G-N_G[v]$.
Hence $H - N_G[v] \in {\rm FS} (\widehat{G - N_G[v]})$.
Since $N_G(v) \subset V_1$,
\begin{itemize}
\item 
$v$ is an isolated vertex in $G[V_2]$;
\item 
$G[S] = (G-N_G[v])[S]$.
\end{itemize}
It then follows that
$G[V_2]$ is the union of $(G-N_G[v])[S]$ and the isolated vertex $v$.
Thus 
we have
$c(G[V_2]) = c((G-N_G[v])[S]) +1 $.
From Lemma~\ref{countlemma}, 
$\mu(H -N_G[v]) = 2^{c(G[V_2])-1} = \mu(H)/ 2$.

Suppose that $N_G[v] = [n-1]$.
%
%
Let $H_0$ be the star subgraph of $G$ with the edge set
$\{\{i, v\} : i \in [n] \setminus \{v\}\}$.
Then $H_0$ belongs to ${\rm FS}_1(\Ghat)$ with $\mu(H_0) = 2$.
We will show that
${\rm FS}_1(\Ghat) =\{H_0\}$.
Suppose that $H \in {\rm FS}_1(\Ghat)$.
Since $N_G[v] = [n-1]$, $V_1=[n] \setminus \{v\}$ and $V_2 = \{v\}$.
Hence $H=H_0$.

\medskip

\noindent
{\bf Proof of (\ref{c4}).} It is enough to show that
 $\varphi: {\rm FS}_1 (\widehat{G-v}) \rightarrow {\rm FS}_1(\Ghat)$,
 $H_1 \mapsto H$ defined below is an injection
such that $\mu(H_1) = \mu(\varphi(H_1)) / 2$.


Let $H_1 \in {\rm FS}_1 (\widehat{G-v})$
and let $[n] \setminus \{v\} = V_1 \sqcup V_2$ where
$n \in V_1$ and $N_G(v) \subset V_1$.
Then $V_2$ is a dominating set of $G-v$.
Let $S=V_2 \cup \{v\}$.
Then 
$$
[n]\setminus \{v\} = \bigcup_{v' \in V_2} N_{G-v }[v'] \subset \bigcup_{v' \in S} N_{G}[v'].
$$
%
Since $v \in S \subset  \bigcup_{v' \in S} N_{G}[v']$, we have 
$$
[n] = V_1 \sqcup S   = \bigcup_{v' \in S} N_{G}[v']
$$
i.e., 
$S$ is a dominating set of $G$.
Hence the bipartite graph $H$ obtained from $H_1$ by adding the vertex $v$ and edges $\{v, v'\}$ ($v' \in N_G(v)$)
is a facet subgraph of $\Ghat$.
Since $N_G(v) \subset V_1$ and $v \notin V_2$,
\begin{itemize}
\item 
$v$ is an isolated vertex in $G[S]$;
\item 
$G[V_2] = (G-v)[V_2]$.
\end{itemize}
Thus $G[S]$ is the union of $(G-v)[V_2]$ and the isolated vertex $v$.
Hence we have $c( G[S] ) = c( (G-v)[V_2] )+1$.
From Lemma~\ref{countlemma}, 
$\mu(H)= 2^{c( (G-v)[V_2] )+1}  = 2 \mu(H_1)$.

\medskip

\noindent
{\bf Proof of (\ref{c5}).} It is enough to show that
$\varphi: {\rm FS}_1 (\widehat{G/v}) \rightarrow {\rm FS}_1(\Ghat)$,
$H_2 \mapsto H$ defined below is an injection 
such that $\mu(H_2) = \mu(\varphi(H_2))/2$.

Let $H_2 \in {\rm FS}_1 (\widehat{G/v})$ 
and let $[n] \setminus \{v\} = V_1 \sqcup V_2$
be the bipartition of the vertex set of $H_2$ with $n \in V_1$ and $N_G(v) \subset V_1$.
Then $V_2$ is a dominating set of $G/v$.
Let $S:= V_2 \cup \{v\}$.
It is trivial that $S \subset \bigcup_{v' \in S} N_G[v']$.
Since $v \in S$, we have $N_G (v)  \subset \bigcup_{v' \in S} N_G[v']$.
Let $v_1 \in V_1 \setminus N_G(v)$.
Since $V_2$ is a dominating set of $G/v$,
there exists $v_2 \in V_2$ such that
$\{v_1,v_2\}$ is an edge of $G/v$.
Since $v_1 \notin N_G(v)$, $\{v_1,v_2\}$ is an edge of $G$.
Hence $V_1 \setminus N_G(v) \subset \bigcup_{v' \in S} N_G[v']$.
Thus 
$$[n] = ((V_1 \setminus N_G(v)) \sqcup N_G(v) ) \sqcup  S = \bigcup_{v' \in S} N_G[v'],$$ i.e.,
$S$ is a dominating set of $G$.
Then
the maximal spanning bipartite subgraph $H$ of $\Ghat$
on the vertex set $V_1 \sqcup S$
is a facet subgraph of $\Ghat$.
Since $N_G(v) \subset V_1$ and $v \notin V_2$,
\begin{itemize}
\item 
$v$ is an isolated vertex in $G[S]$;
\item 
$G[V_2] = (G/v)[V_2]$.
\end{itemize}
It then follows that
$G[S]$ is the union of $(G/v)[V_2]$ and the isolated vertex $v$.
Thus we have $c( G[S] ) = c(  (G/v) [V_2])+1$.
From Lemma~\ref{countlemma}, 
$\mu(H)= 2^{c(  (G/v) [V_2])+1} = 2 \mu (H_2) $.

\bigskip

Using (\ref{c1}) -- (\ref{c5}),
we will show 
(\ref{keyF}).
%
%
%
From (\ref{c1}), (\ref{c2}) and (\ref{c3}), 
we have
\begin{align*}
N(\Pc_\Ghat) 
&= \sum_{H \in {\rm FS}_0(\Ghat)} \mu(H)
+ \sum_{H \in {\rm FS}_1(\Ghat)} \mu(H)
+ \sum_{H \in {\rm FS}_2(\Ghat)} \mu(H)\\
&\le 
\sum_{H \in {\rm FS}_2(\widehat{G-v})} \mu(H)
+ \sum_{H \in {\rm FS}(\widehat{G-N_G[v]})} 2 \mu(H)
+ \sum_{H \in {\rm FS}_2(\widehat{G/v})} \mu(H)\\
&\le 
\sum_{H \in {\rm FS}(\widehat{G-v})} \mu(H)
+ \sum_{H \in {\rm FS}(\widehat{G-N_G[v]})} 2 \mu(H)
+ \sum_{H \in {\rm FS}(\widehat{G/v})} \mu(H)\\
&= 
    N(\Pc_{\widehat{G-v}})+ 2N(\Pc_{\widehat{G-N_G[v]}})+ N(\Pc_{\widehat{G/v}}).  
\end{align*}
%
Moreover,
from (\ref{c1}), (\ref{c2}), (\ref{c4}) and (\ref{c5}), 
we have
\begin{align*}
N(\Pc_{\widehat{G-v}})  + N(\Pc_{\widehat{G/v}})
&= 
\sum_{H \in {\rm FS}_1(\widehat{G-v})} \mu(H)
+ \sum_{H \in {\rm FS}_2(\widehat{G-v})} \mu(H)
+\sum_{H \in {\rm FS}_1(\widehat{G/v})} \mu(H)
+ \sum_{H \in {\rm FS}_2(\widehat{G/v})} \mu(H)
\\
&\le
\sum_{H \in {\rm FS}_1(\widehat{G})} \frac{1}{2} \mu(H)
+ \sum_{H \in {\rm FS}_0(\widehat{G})} \mu(H)
+ \sum_{H \in {\rm FS}_1(\widehat{G})} \frac{1}{2} \mu(H)
+\sum_{H \in {\rm FS}_2(\widehat{G})} \mu(H)
\\
&= 
    N(\Pc_{\widehat{G}}).
\end{align*}

Finally , we will show (\ref{keyF2}).
%
%
From (\ref{c1}), (\ref{c2}) and (\ref{c3}), 
we have
\begin{align*}
N(\Pc_{\widehat{G}}) &= 
\sum_{H \in {\rm FS}_0(\widehat{G})} \mu(H)
+\sum_{H \in {\rm FS}_2(\widehat{G})} \mu(H) +2\\
&=
\sum_{H \in {\rm FS}_2(\widehat{G-v})} \mu(H)
+\sum_{H \in {\rm FS}_2(\widehat{G/v})} \mu(H)+2
\\
&= 
N(\Pc_{\widehat{G-v}})  + N(\Pc_{\widehat{G/v}})+2,
\end{align*}
as desired.
\end{proof}

\begin{Corollary}
    Let $G$ be a graph with $n-1 \ge 2$ vertices.
    Then 
    \[
N(\Pc_{\widehat{\widehat{G\;}}})
    =
    N(\Pc_\Ghat) + 2^n.
    \]
\end{Corollary}

\begin{proof}
Note that $\widehat{G}$ has a vertex $v$ of degree $n-1$.
Then $\widehat{G} -v = G$ and $\widehat{G}/v = K_{n-1}$.
From Proposition~\ref{insp} (\ref{keyF2}),
we have 
\[N(\Pc_{\widehat{\widehat{G\;}}})
=N(\Pc_\Ghat) + N(\Pc_{\widehat{K_{n-1}}}) + 2
=N(\Pc_\Ghat) + (2^n -2) + 2
=
N(\Pc_\Ghat) + 2^n.\]    
\end{proof}

We are in a position to prove Theorem~\ref{mainT}.

\begin{proof}[Proof of Theorem \ref{mainT}]
Proof is by induction on $n$ ($\ge 2$).
If $n=2$, then $\Ghat = K_2 $ and hence
$$2^1 = N(\Pc_\Ghat)  < 14 \cdot 6^{-1}.$$
Thus the assertion holds.
Suppose that $n>2$ and the assertion is true for the graphs with less number of vertices.

\medskip

\noindent
{\bf Case 1} ($G$ has no vertices of degree $\ge 2$){\bf .}
Then $G$ is a disjoint union of edges $e_1,\dots,e_t$
and isolated vertices $v_1,\dots,v_{n-2t -1}$.
The suspension $\Ghat$ of $G$ is a $1$-sum of $t$ triangles
together with $n-2t-1$ edges.
Since $N(\Pc_{K_2})=2$ and $N(\Pc_{K_3})=6$, we have 
$$N(\Pc_{\Ghat}) = 2^{n-2t-1} \cdot 6^t = \left( \frac{2}{3} \right)^{\frac{n-2t-1}{2}} 6^{\frac{n-1}{2}} \le 6^{\frac{n-1}{2}}.$$
The equality holds if and only if $n-2t -1=0$, that is, 
$\Ghat$ is a $1$-sum of $t$ triangles.
Note that $n$ is odd if $n-2t -1=0$.
Suppose that $n$ ($\ge 4$) is even.
Then $n-2t -1 \ge 1$.
Hence 
$$
\left( \frac{2}{3} \right)^{\frac{n-2t-1}{2}} 6^{\frac{n-1}{2}}
\le 
\left( \frac{2}{3} \right)^{\frac{1}{2}} 6^{\frac{n-1}{2}}
=
12 \cdot 6^{\frac{n}{2}-2}
<
14 \cdot 6^{\frac{n}{2}-2}.
$$

On the other hand, 
$$N(\Pc_{\Ghat}) = 2^{n-2t-1} \cdot 6^t = \left( \frac{3}{2} \right)^t
 2^{n-1}
 \ge 
2^{n-1}
 ,$$
and equality holds if and only if $t=0$, that is,
$G$ is an empty graph.

\medskip

\noindent
{\bf Case 2} ($G$ has a vertex $v$ of degree $\ge 2$){\bf .}
Then $n \ge 4$.
    Since $\deg(v) \ge 2$, $G/v$ is not empty.
By the hypothesis of induction, 
$$
\begin{array}{cclccc}
2^{n-2} &\le&  N(\Pc_{\widehat{G-v}})& \le & 6^{\frac{n-2}{2}},&\\
2^{n-2} &<&  N(\Pc_{\widehat{G/v}}) &\le & 6^{\frac{n-2}{2}},&\\
&  &  N(\Pc_{\widehat{G-N_G[v]}}) &\le & 6^{\frac{n-4}{2}} &  (\mbox{if } N_G[v] \ne [n-1]).
\end{array}
$$

\medskip

\noindent
\textbf{Case 2.1} ($N_G[v] = [n-1]$){\bf .}
From Proposition~\ref{insp} (\ref{keyF2}), 
we have 
\begin{align*}
         N(\Pc_\Ghat) &=   N(\Pc_{\widehat{G-v}})+  N(\Pc_{\widehat{G/v}}) +2
 >
 2^{n-2} + 2^{n-2} + 2 >2^{n-1}, \\
  N(\Pc_\Ghat) &=   N(\Pc_{\widehat{G-v}})+  N(\Pc_{\widehat{G/v}}) +2
 \le 6^{\frac{n-2}{2}} + 6^{\frac{n-2}{2}} +2
 = 12 \cdot 6^{\frac{n}{2}-2} +2
 \le 14 \cdot 6^{\frac{n}{2}-2} <  6^{\frac{n-1}{2}}.
\end{align*}
In addition, 
$N(\Pc_\Ghat) =  14 \cdot 6^{\frac{n}{2}-2} $
if and only if $n=4$ and 
both $G-v$ and $G/v$ are $K_2$
if and only if $G$ is a triangle.

\medskip

\noindent
\textbf{Case 2.2} ($N_G[v] \ne [n-1]$){\bf .}
From Proposition~\ref{insp} (\ref{keyF}), we have 
\[
N(\Pc_\Ghat) \ge N(\Pc_{\widehat{G-v}})+ N(\Pc_{\widehat{G/v}})
>
2^{n-2}   +2^{n-2} =2^{n-1},
\]
and
\begin{eqnarray}
N(\Pc_\Ghat) &\le& 
    N(\Pc_{\widehat{G-v}})+ 2N(\Pc_{\widehat{G-N_G[v]}})+ N(\Pc_{\widehat{G/v}})  \label{one}\\
    &\le & 6^{\frac{n-2}{2}}  + 2 \cdot 6^{\frac{n-4}{2}}+  6^{\frac{n-2}{2}}\\
    &=& \frac{7}{3} \cdot 6^{\frac{n-2}{2}} = 14 \cdot 6^{\frac{n}{2}-2} <  6^{\frac{n-1}{2}}. \label{three}
\end{eqnarray}

Suppose that $n$ is even.
We will show that $N(\Pc_\Ghat)= 14 \cdot 6^{\frac{n}{2}-2}$
if and only if $G$ is a disjoint union of several edges with the triangle.

\medskip

\noindent
(If) 
If $G$ is a disjoint union of several edges with the triangle,
then $\widehat{G}$ is a disjoint union of several triangles with $K_4$.
Since $N(\Pc_{K_3}) = 6$ and $N(\Pc_{K_4})=14$, $N(\Pc_\Ghat)= 14 \cdot 6^{\frac{n}{2}-2}$.

\medskip

\noindent
(Only if)
Suppose that $N(\Pc_\Ghat)= 14 \cdot 6^{\frac{n}{2}-2}$.
From (\ref{one}) -- (\ref{three}) above, by the hypothesis of induction,
each of $G-v$, $G/v$ and $G - N_G[v]$ 
is a disjoint union of several edges, and the number of vertices of $G-N_G[v]$ is $n-4$.
Then $\deg(v) =2$.
Let $N_G(v) = \{v_1, v_2\}$.
Since $G/v$ is a disjoint union of several edges
and since $\{v_1,v_2\}$ is an edge of $G/v$,
$N_{G/v}(v_1) = \{v_2\}$ and 
$N_{G/v}(v_2) = \{v_1\}$.
Since $G-v$ is a disjoint union of several edges, $\{v_1, v_2\}$ is an edge of $G$.
In addition, since $G - N_G[v]$
is a disjoint union of several edges,
$G$ is a disjoint union of several edges with the triangle 
$(v,v_1,v_2)$.
\end{proof}

\begin{Remark}
Given a graph $G$ on the vertex set $[n]$,
let $Q_{ij}(G)$ denote the number of subset 
$S \subset [n]$ with $i=|S|$ and $j = c(G[S])$.
Then the polynomial 
\[Q(G;x,y) = \sum_{i=0}^n \sum_{j=0}^n Q_{ij}(G) x^i y^j
\]
is called the {\it subgraph component polynomial} of $G$.
From Lemma~\ref{countlemma}, it follows that
$Q(G;1,2)$ gives an upper bound of $N(\Pc_\Ghat)$.
Although it seems to be difficult to apply the theory of 
subgraph component polynomials to our problem directly, the idea of the proof of Proposition~\ref{insp} is inspired by 
\cite[Theorem 13]{TAM}.

\end{Remark}

\section{Join graphs}
\label{sec:four}

In the present section, we extend Theorem~\ref{mainT} to join graphs
by giving a proof of Theorem~\ref{JOIN theorem}.

\begin{Lemma} \label{join to suspension}
    Let $G_1 =(V, E)$ and $G_2=(V', E')$ be graphs with $V \cap V' = \emptyset$, $|V|=n_1$, and $|V'|=n_2$.
    For each $i=1,2$, let $m_i$ be the number of connected components of $G_i$.
    Then we have
    \[
    N(\Pc_{G_1+ G_2}) 
    \ \ 
    \le \ \ 
    N(\Pc_{\widehat{G_1}})+ N(\Pc_{\widehat{G_2}}) +2^{m_1} +2^{m_2} -2 +4(2^{n_1-1} -1)(2^{n_2-1} -1).
    \]
\end{Lemma}

\begin{proof}
We define a partition 
${\rm FS} (G_1 + G_2) = {\rm FS}_1 \sqcup {\rm FS}_2 \sqcup {\rm FS}_3  \sqcup {\rm FS}_4$, where 
\begin{eqnarray*}
{\rm FS}_1  &:=& \left\{ H \in {\rm FS}(G_1 + G_2) : 
\begin{array}{c} \mbox{the bipartition of } H \mbox{ is } V_1 \sqcup V_2, \mbox{ where }\\ 
V \cap V_1\ne \emptyset, V \cap V_2 \ne \emptyset \mbox{ and } V' \subset V_1\end{array} 
\right\},\\
{\rm FS}_2  &:=& \left\{ H \in {\rm FS}(G_1 + G_2) : 
\begin{array}{c} \mbox{the bipartition of } H \mbox{ is } V_1 \sqcup V_2, \mbox{ where }\\ 
V' \cap V_1 \ne \emptyset, V' \cap V_2 \ne \emptyset  \mbox{ and } V \subset V_1\end{array} 
\right\},\\
{\rm FS}_3  &:=& \left\{ H \in {\rm FS}(G_1 + G_2) : 
 \mbox{the bipartition of } H \mbox{ is } V \sqcup V'
\right\},\\
{\rm FS}_4  &:=& \left\{ H \in {\rm FS}(G_1 + G_2) : 
\begin{array}{c} \mbox{the bipartition of } H \mbox{ is } V_1 \sqcup V_2, \mbox{ where }\\ 
V \cap V_1\ne \emptyset,  V \cap V_2\ne \emptyset, V' \cap V_1\ne \emptyset, V' \cap V_2 \ne \emptyset\end{array} 
\right\}.
\end{eqnarray*}

\medskip

\noindent
{\bf Claim 1.} There is an injection $\varphi: {\rm FS}_1 \rightarrow {\rm FS}(\widehat{G_1})$
such that $\mu(H) = \mu(\varphi(H))$.

Let $H \in {\rm FS}_1$.
Then $V \cap V_2$ is a dominating set of $G_1$.
Hence the graph $H'$ obtained from $H$ by contracting the vertices in $V'$ to one vertex
is a facet subgraph of $\widehat{G_1}$.
Since $(G_1+G_2)[V_1]$ is connected, we have $\mu(H) = \mu(H')
=2^{c( G_1 [V_2])}$.

\medskip

\noindent
{\bf Claim 2.}  There is an injection $\varphi: {\rm FS}_2 \rightarrow {\rm FS}(\widehat{G_2})$
such that $\mu(H) = \mu(\varphi(H))$.

It follows from the same argument as in Claim 1.

\medskip

\noindent
{\bf Claim 3.} $ {\rm FS}_3=\{H_0\}$ where $\mu(H_0)=2^{m_1} +2^{m_2} -2$.

Let $H'$ denote the graph obtained from $H$ by 
contracting each connected component of $G_1[V]$ and that of $G_2[V']$ to one vertex.
From Proposition~\ref{facet defining}, $\mu(H_0) = N(\Pc_{H'})$. 
Since $H'$ is a complete bipartite graph with partition $V_1' \sqcup V_2'$,
where $|V_1'|=m_1$ and $|V_2'|=m_2$, it follows from equation (\ref{2bu graph}) in Introduction that $\mu(H_0)=2^{m_1} +2^{m_2} -2$.

\medskip

\noindent
{\bf Claim 4.} $|{\rm FS}_4| \le 2(2^{n_1-1} -1)(2^{n_2-1} -1)$ and $\mu(H)=2$ for each $H \in {\rm FS}_4$.

The number of facet subgraphs $H \in {\rm FS}_4$ is at most $2(2^{n_1-1} -1)(2^{n_2-1} -1)$ by considering the possibility of $V_1$ and $V_2$.
If $H \in {\rm FS}_4$, then both $(G_1+G_2)[V_1]$ and $(G_1+G_2)[V_2]$ are connected, and hence $\mu(H)=2$ from Proposition~\ref{facet defining}.

\medskip

\noindent
From Claims 1, 2, 3, and 4 we have 
\[    N(\Pc_{G_1+ G_2}) 
    \ \ 
    \le \ \ 
    N(\Pc_{\widehat{G_1}})+ N(\Pc_{\widehat{G_2}}) +2^{m_1} +2^{m_2} -2 +4(2^{n_1-1} -1)(2^{n_2-1} -1),
\]
as desired.
\end{proof}

We now prove the main theorem of the present paper.

\begin{proof}[Proof of Theorem~\ref{JOIN theorem}]
From Proposition~\ref{small}, we may assume that $n \ge 5$.
Let $n_1=|V|$ and $n_2 = |V'|$.
From Theorem~\ref{mainT},
we may assume that $G_1 + G_2$ has no vertices of degree $n-1$.
In addition, if both $G_1$ and $G_2$ are empty, then 
$G_1 + G_2$ is a complete bipartite graph and hence satisfies the assertion.
Thus we may assume that
\begin{itemize}
    \item[(i)] each $G_i$ has no vertices of degree $n_i -1$,
    \item[(ii)] $n_1 \ge n_2 \ge 2$, and $n\ge 5$,
    \item[(iii)] either $G_1$ or $G_2$ has at least one edge. 
\end{itemize}

First, we will show 
$N(\Pc_{G_1 + G_2}) > 3  \cdot 2^{\frac{n-1}{2}}-2 \ (> 2^{\frac{n}{2}+1}-2)$.
Let $ {\rm FS}_3=\{H_0\}$ and ${\rm FS}_4$ denote the sets defined in the proof of Lemma~\ref{join to suspension}. 
Let $\{i,j\}$ be an edge of $G_1$.
Then a maximal spanning bipartite subgraph of $G_1 + G_2$ with partition 
$V_1 \sqcup V_2$ where $i \in V_1$, $j \in V_2$,
$V_1 \cap V' \ne \emptyset$ and $V_2 \cap V' \ne \emptyset$
belongs to ${\rm FS}_4$.
The number of such partitions equals to $2^{n_1-2} (2^{n_2}-2)=2^{n-2} - 2^{n_1-1}$.
Hence $|{\rm FS}_4| \ge 2^{n-2} - 2^{n_1-1}$.
Similarly, if $G_2$ has an edge,
then $|{\rm FS}_4| \ge 2^{n-2} - 2^{n_2-1}$.
Since $n-2 \ge n_1 \ge n_2$, we have $|{\rm FS}_4| \ge 2^{n-2} - 2^{n-3} = 2^{n-3}$.
Then
\[
N(\Pc_{G_1 + G_2})  \ge 2 \cdot |{\rm FS}_4| +\mu (H_0) 
\ge 2^{n-2} + 2^{m_1} + 2^{m_2} -2 \ (> 2^{n-2}).
\]
If $n=5$, then $(n_1,n_2) = (3,2)$ and $G_2$ is an empty graph with 2 vertices.
Since $m_2 =2$ and $m_1 \ge 1$,
\[
(2^{n-2} + 2^{m_1} + 2^{m_2} -2) -( 3  \cdot 2^{\frac{n-1}{2}}-2) \ge 2 >0.
\]
If $n=6$, then 
$2^{n-2} -( 3  \cdot 2^{\frac{n-1}{2}}-2) = 6(3 - 2 \sqrt{2}) >0.$
If $n\ge 7$, then 
\[
2^{n-2} -( 3  \cdot 2^{\frac{n-1}{2}}-2) = 2^{\frac{n-1}{2}} \left(2^{\frac{n-3}{2}} -3 \right)+2
>0.
\]
Thus we have $N(\Pc_{G_1 + G_2}) > 3  \cdot 2^{\frac{n-1}{2}}-2 $.

Finally, we will show $N(\Pc_{G_1 + G_2}) <  14 \cdot 6^{\frac{n}{2}-2} \ 
(< 6^{\frac{n-1}{2}})$.

\medskip

\noindent
{\bf Case 1}
($n_2=2$){\bf .}
From (i) above, $G_2$ is an empty graph with 2 vertices
and hence $N(\Pc_{\widehat{G_2}}) = 4$.
From (iii), $G_1$ has at least one edge.
In particular, the number of connected components of $G_1$ is $m_1 < n_1 =n-2$.
From 
Lemma~\ref{join to suspension},
\begin{align*}
 N(\Pc_{G_1+ G_2}) 
    &\le 
    N(\Pc_{\widehat{G_1}})+ N(\Pc_{\widehat{G_2}}) +2^{m_1} +2^{m_2} -2 +4(2^{n_1-1} -1)(2^{n_2-1} -1)\\
    &\le 
   N(\Pc_{\widehat{G_1}}) + 4  +2^{n-3} +2^{2} -2 +4(2^{n-3} -1)(2^{2-1} -1)\\
   &=N(\Pc_{\widehat{G_1}})
   + 5 \cdot 2^{n-3} +    2.
\end{align*}
If $n=5$, then 
$(n_1,n_2) = (3,2)$ and $G_2$ is an empty graph with 2 vertices.
From (i) and (iii) above, $G_1$ has exactly one edge.
Thus $N(\Pc_{\widehat{G_1}})=12$, and hence
$N(\Pc_{\widehat{G_1}})
   + 5 \cdot 2^{n-3} +    2  = 34 < 14 \sqrt{6}$.
Suppose that $n  \ge 6$.
From Theorem~\ref{mainT},
$$
14 \cdot 6^{\frac{n}{2}-2} - (6^{\frac{n}{2}-1} + 5 \cdot 2^{n-3} +    2)
=
48 \cdot 6^{\frac{n}{2}-3} -  40 \cdot 4^{\frac{n}{2}-3} - 2 >0.
$$
Thus we have $N(\Pc_{G_1+ G_2})  < 14 \cdot 6^{\frac{n}{2}-2}$.

\medskip

\noindent
{\bf Case 2}
($n_2 \ge 3$){\bf .}
Then $n\ge 6$.
From Theorem~\ref{mainT} and Lemma~\ref{join to suspension},
\begin{align*}
 N(\Pc_{G_1+ G_2}) 
    &\le 
    N(\Pc_{\widehat{G_1}})+ N(\Pc_{\widehat{G_2}}) +2^{m_1} +2^{m_2} -2 +4(2^{n_1-1} -1)(2^{n_2-1} -1)\\
    &\le 
   6^{\frac{n_1}{2}} + 6^{\frac{n_2}{2}}  +2^{n_1} +2^{n_2} -2 +4(2^{n_1-1} -1)(2^{n_2-1} -1)\\
   &
   = 6^{\frac{n_1}{2}} + 6^{\frac{n_2}{2}} + 2^{n_1+n_2} -2^{n_1} -2^{n_2} + 2\\
   &
   \le 2 \cdot 6^{\frac{n-3}{2}} + 2^{n} -14.
\end{align*}
If $n =6$, then $14 \cdot 6^{\frac{n}{2}-2}
-
(2 \cdot 6^{\frac{n-3}{2}} + 2^{n} -14)
= 34 -12\sqrt{6} > 0$.
If $n \ge 7$, then we have
\[
14 \cdot 6^{\frac{n}{2}-2}
-
(2 \cdot 6^{\frac{n-3}{2}} + 2^{n} -14)
=
(84 \sqrt{6} - 72)  \cdot 6^{\frac{n-7}{2}}
-
128 \cdot 4^{\frac{n-7}{2}} +14
>0.
\]
(Here, $84 \sqrt{6} - 72 \risingdotseq 133.76$.)
Thus we have $N(\Pc_{G_1+ G_2})  < 14 \cdot 6^{\frac{n}{2}-2}$.
\end{proof}

In the present paper, 
we proved that Conjecture \ref{yoso} is true
for any graph that is the join of two graphs.
The proofs depend on the structure of such graphs, i.e, 
there exists a vertex with relatively large degree, and hence
$\mu(H)$ is relatively easy to compute for each facet subgraph $H$.

From Theorem~\ref{JOIN theorem}, Propositions~\ref{tree facet} and \ref{block prop}, in order to study Conjecture \ref{yoso},
it is enough to discuss $2$-connected non-bipartite graphs whose complement is connected.

\section*{Acknowledgement}
The authors are grateful to anonymous referees for their careful reading and helpful comments.
This work was supported by 
Grant-in-Aid for JSPS Fellows 23KJ2117,
and
JSPS KAKENHI 24K00534.

\end{document}